\def\ifundefined#1#2{\expandafter\ifx\csname#1\endcsname\relax\input #2\fi}
\input amssym.def

\def\pmb#1{\setbox0=\hbox{$#1$}       
     \kern-.025em\copy0\kern-\wd0
     \kern.05em\copy0\kern-\wd0
     \kern-.025em\box0}




\def\Cross{\bigm| \kern-5.5pt \not \ \, }
\def\cross{\mid \kern-5.0pt \not \ \, }             
\def\notto{\hbox{$~\rightarrow~\kern-1.5em\hbox{/}\ \ $}}

\hyphenation{math-ema-ticians}
\hyphenation{pa-ra-meters}
\hyphenation{pa-ra-meter}
\hyphenation{lem-ma}
\hyphenation{lem-mas}
\hyphenation{to-po-logy}
\hyphenation{to-po-logies}
\hyphenation{homo-logy}
\hyphenation{homo-mor-phy}

\def\nSigma{\Sigma \kern-8.3pt \bigm|\,}

\font\teneufm=eufm10
\font\eighteufm=eufm8
\font\fiveeufm=eufm5

\newfam\eufam
\textfont\eufam=\teneufm
\scriptfont\eufam=\eighteufm
\scriptscriptfont\eufam=\fiveeufm

\def\boxit#1{\vbox{\hrule\hbox{\vrule\kern2.0pt
       \vbox{\kern2.0pt#1\kern2.0pt}\kern2.0pt\vrule}\hrule}}
\def\vlra#1{\hbox{\kern-1pt
       \hbox{\raise2.38pt\hbox{\vbox{\hrule width#1 height0.26pt}}}
       \kern-4.0pt$\rightarrow$}}

\def\vlla#1{\hbox{$\leftarrow$\kern-1.0pt
       \hbox{\raise2.38pt\hbox{\vbox{\hrule width#1 height0.26pt}}}}}

\def\vlda#1{\hbox{$\leftarrow$\kern-1.0pt
       \hbox{\raise2.38pt\hbox{\vbox{\hrule width#1 height0.26pt}}}
       \kern-4.0pt$\rightarrow$}}

\def\longra#1#2#3{\,\lower3pt\hbox{${\buildrel\mathop{#2}
\over{{\vlra{#1}}\atop{#3}}}$}\,}

\def\longla#1#2#3{\,\lower3pt\hbox{${\buildrel\mathop{#2}
\over{{\vlla{#1}}\atop{#3}}}$}\,}

\def\longda#1#2#3{\,\lower3pt\hbox{${\buildrel\mathop{#2}
\over{{\vlda{#1}}\atop{#3}}}$}\,}

\def\overrightharpoonup#1{\vbox{\ialign{##\crcr
	$\rightharpoonup$\crcr\noalign{\kern-1pt\nointerlineskip}
	$\hfil\displaystyle{#1}\hfil$\crcr}}}

\catcode`@=11
\def\@@dalembert#1#2{\setbox0\hbox{$#1\rm I$}
  \vrule height.90\ht0 depth.1\ht0 width.04\ht0
  \rlap{\vrule height.90\ht0 depth-.86\ht0 width.8\ht0}
  \vrule height0\ht0 depth.1\ht0 width.8\ht0
  \vrule height.9\ht0 depth.1\ht0 width.1\ht0 }
\def\dalembert{\mathord{\mkern2mu\mathpalette\@@dalembert{}\mkern2mu}}

\def\@@varcirc#1#2{\mathord{\lower#1ex\hbox{\m@th${#2\mathchar\hex0017 }$}}}
\def\varcirc{\mathchoice
  {\@@varcirc{.91}\displaystyle}{\@@varcirc{.91}\textstyle}
{\@@varcirc{.45}\scriptscriptstyle}}
\catcode`@=12

\font\tensf=cmss10 \font\sevensf=cmss8 at 7pt
\newfam\sffam
\textfont\sffam=\tensf\scriptfont\sffam=\sevensf

\input amssym.def
\input amssym
\magnification=1200

\font\bigslll=cmsl10 scaled\magstep4
\tolerance=500
\overfullrule=0pt
\null
\centerline{\bigslll Generalized Serre--Tate Ordinary Theory}
\bigskip\bigskip\bigskip\medskip
\centerline{Adrian Vasiu${}^*$ $\vfootnote{*} {Department of Mathematical Sciences, Binghamton University, Binghamton, P. O. Box 6000, NY 13902-6000, U.S.A. {} {} e-mail: adrian\hbox{@}math.binghamton.edu}$  }
\bigskip\bigskip\medskip
\centerline{final version, April 30, 2012, to be published by International Press, Inc.} 
\bigskip\bigskip\noindent
{\bf Abstract}. We study a generalization of Serre--Tate theory of ordinary abelian varieties and their deformation spaces. This generalization deals with abelian varieties equipped with additional structures. The additional structures can be not only an action of a semisimple algebra and a polarization, but more generally the data given by some ``crystalline Hodge cycles'' (a $p$-adic version of a Hodge cycle in the sense of motives). Compared to Serre--Tate ordinary theory, new phenomena appear in the generalized context. 

We give an application of the generalized theory to the existence of ``good'' integral models of those Shimura varieties whose adjoints are products of simple, adjoint Shimura varieties of $D_l^{\bf H}$ type with $l\ge 4$.  
\bigskip\medskip\medskip\noindent
{\bf Key words}: $p$-divisible groups and objects, $F$-crystals, Newton polygons, ordinariness, canonical lifts, reductive group schemes, formal Lie groups, abelian varieties, Shimura varieties, integral models, Hodge cycles, complex multiplication, and deformation spaces.
\bigskip\medskip\medskip\noindent
{\bf MSC 2000}: 11G10, 11G15, 11G18, 14F20, 14F30, 14F40, 14G35, 14K10, 14K22, 14L05, 14L15, and 17B45. 

\bigskip\bigskip\bigskip\bigskip\bigskip\bigskip\bigskip\medskip\medskip\medskip\noindent
\centerline{All copyrights reserved to International Press, Inc.}

\bigskip\noindent
Final version 196 pages (including contents and index) to be published by International Press, Inc. Until publication the pdf file is available at http://www.math.binghamton.edu/adrian/" 
\end